\documentclass[11pt]{article}
\usepackage[francais]{babel}
\usepackage[latin1]{inputenc}
\usepackage[T1]{fontenc}
\usepackage{a4,amsmath,amsfonts}
\title{Géométrie des surfaces algébriques et points entiers}
\author{Pascal Autissier}
\begin{document}

\maketitle

\newcommand{\D}{\displaystyle}

{\bf Abstract:} Let $X$ be a projective normal surface over a number field $K$.
Let $H$ be the sum of four properly intersecting ample effective divisors on
$X$. We show that any set of $S$-integral points in $X-H$ is not Zariski
dense.\\


\section{Introduction}

On s'intéresse ici aux solutions à coordonnées quasi-entières de systèmes
d'équations polynomiales à coefficients dans un corps de nombres, dans l'esprit
de la conjecture de Lang et Vojta ({\it cf} conjecture 4.2 de \cite{Lang}
p. 223).

Plus précisément, soient $K$ un corps de nombres et $S$ un ensemble fini de
places de $K$. On montre le résultat suivant:\\

{\bf Théorème 1.1:} {\it Soit $X$ une surface normale projective sur $K$.
Soient $D_1;D_2;D_3;D_4$ quatre diviseurs effectifs amples sur $X$ qui se
coupent proprement. Posons $Y=X-D_1\cup D_2\cup D_3\cup D_4$. Soit
${\cal E}\subset Y(K)$ un ensemble $S$-entier sur $Y$. Alors ${\cal E}$ n'est
pas Zariski-dense dans $Y$.}\\

Cet énoncé était connu de Vojta pour $X$ lisse vérifiant $\rho\le g+1$, où
$\rho$ désigne le nombre de Picard de $X_{\overline{K}}$ et
$g=h^1(X;{\cal O}_X)$ ({\it cf} corollaire 0.3 de \cite{Voj2}). En fait, Vojta
a besoin de $\rho+3-g$ diviseurs au lieu de 4. L'intérêt de notre résultat
réside donc dans l'uniformité en le nombre de diviseurs à considérer.\\

Remarquons que le théorème 1.1 s'inscrit bien dans le cadre de la conjecture
de Lang et Vojta, puisque si $X$ est lisse sur $K$ de diviseur canonique
${\cal K}_X$, alors ${\cal K}_X+D_1+D_2+D_3+D_4$ est ample sur $X$ ({\it cf}
exemple 1.5.35 de \cite{Laza} p. 87).\\

La démonstration repose sur une légère extension d'un théorème de Corvaja et
Zannier \cite{CoZa} ({\it cf} théorème 3.2), qui donne des conditions
géométriques de non-Zariski-densité des points $S$-entiers, et sur un bon choix
de multiplicités associées aux diviseurs $D_i$ ({\it cf} proposition 2.3).

On utilise en particulier le théorème du sous-espace de Schmidt ({\it cf}
\cite{Schm} \S VI) et Schlickewei ({\it cf} \cite{Schl}).\\

Je remercie Antoine Chambert-Loir pour l'inspiration qu'il m'a procurée.\\

\section{Géométrie}

Soit $K$ un corps de caractéristique nulle.\\

{\bf Conventions:} On appelle variété sur $K$ tout schéma intègre,
quasi-projectif et géométriquement irréductible sur $K$. Une surface sur $K$
est une variété sur $K$ de dimension 2. Le mot ``diviseur'' sous-entend
``diviseur de Cartier''.\\

Soit $X$ une variété projective sur $K$ de dimension $d\ge1$. Lorsque
$L_1;\cdots;L_d$ sont des diviseurs sur $X$, on désigne par
$\bigl<L_1\cdots L_d\bigr>$ leur nombre d'intersection.\\

Soient $L$ un diviseur ample sur $X$ et $E$ un diviseur effectif non nul sur
$X$. La formule de Hirzebruch-Riemann-Roch donne l'estimation
$\D h^0(X;nL)=\frac{\bigl<L^d\bigr>}{d!}n^d+O(n^{d-1})$. Ceci motive la
définition suivante:\\

Pour tout entier $n\ge1$, posons d'abord $\D S_n=\sum_{k\ge1}h^0(X;nL-kE)$;
remarquons que cette somme est finie puisque $h^0(X;nL-kE)=0$ si
$k>\bigl<L^d\bigr>n/\bigl<L^{d-1}E\bigr>$.\\

{\bf Définition:} On pose $\D\nu(L;E)=\liminf_{n\rightarrow+\infty}\frac{S_n}{h^0(X;nL)n}=\liminf_{n\rightarrow+\infty}\frac{d!S_n}{\bigl<L^d\bigr>n^{d+1}}$.\\

{\bf Exemple:} Si $X$ est une courbe, alors on peut aisément expliciter cette
constante; on trouve $\D\nu(L;E)=\frac{\bigl<L\bigr>}{2\bigl<E\bigr>}$.\\

On aura besoin dans la suite de minorer $\nu(L;E)$. Commençons par une variante
des ``inégalités de Morse holomorphes'' ({\it cf} \cite{Dema} \S 12):\\

{\bf Lemme 2.1:} {\it Soit $X$ une surface projective sur $K$. Soient $L$ et
$M$ des diviseurs amples sur $X$. Posons
$\alpha=\bigl<LM\bigr>/\bigl<M^2\bigr>$. Soient $n$ et $k$ des entiers
vérifiant $1\le k\le\alpha n$. On a alors la minoration
$$h^0(X;nL-kM)\ge\bigl<L^2\bigr>\frac{n^2}{2}-\bigl<LM\bigr>nk+\bigl<M^2\bigr>\frac{k^2}{2}-O(n)\quad,$$
où le $O$ ne dépend que de $(K;X;L;M)$.}\\

{\it Démonstration:} On choisit un entier $b\ge1$ tel que $bM$ soit très ample.
D'après le théorème de Bertini, il existe $s\in\Gamma(X;bM)-\{0\}$ tel que
$C={\rm div}(s)$ soit géométriquement irréductible sur $K$.\\

Soit $i$ un entier tel que $0\le i\le\alpha n$. On a la
suite exacte de ${\cal O}_X$-modules suivante:
$$0\rightarrow{\cal O}_X(nL-(i+b)M)\rightarrow{\cal O}_X(nL-iM)\rightarrow{\cal O}_X(nL-iM)_{|C}\rightarrow0\quad.$$

On en déduit une suite exacte en cohomologie qui donne l'inégalité
$$h^0(X;nL-(i+b)M)\ge h^0(X;nL-iM)-h^0(C;(nL-iM)_{|C})\quad.$$

En utilisant la majoration $h^0(C;L'_{|C})\le\bigl<L'C\bigr>+1$ valable pour
tout diviseur $L'$ tel que $\bigl<L'C\bigr>\ge0$ ({\it cf} proposition 3 (3) de
\cite{Ful1} p. 192), on obtient
$$h^0(X;nL-(i+b)M)\ge h^0(X;nL-iM)-\bigl<LM\bigr>bn+\bigl<M^2\bigr>bi-1\quad.$$

Maintenant, écrivons $k$ sous la forme $k=bq+r$ avec $q\ge0$ et $0\le r<b$. En
sommant l'inégalité précédente, on trouve
$$\begin{array}{rcl}
\D h^0(X;nL-kM)&\ge&\D h^0(X;nL-rM)-\sum_{j=0}^{q-1}\Bigl[\bigl<LM\bigr>bn-\bigr<M^2\bigr>b(bj+r)+1\Bigr]\\
&=&\D\bigl<L^2\bigr>\frac{n^2}{2}-\bigl<LM\bigr>nk+\bigl<M^2\bigr>\frac{k^2}{2}-O(n)\\
\end{array}$$
(l'asymptotique $h^0(X;nL-rM)=\bigl<L^2\bigr>n^2/2+O(n)$ est fournie par
Hirzebruch-Riemann-Roch). D'où le résultat. $\square$\\

{\bf Remarque:} La démonstration donne en fait une minoration de
$h^0(nL-kM)-h^1(nL-kM)$.\\

{\bf Corollaire 2.2:} {\it Soit $X$ une surface projective sur $K$. Soient $L$
un diviseur ample sur $X$ et $E$ un diviseur effectif ample sur $X$. On a alors
$$\nu(L;E)\ge\frac{\bigl<L^2\bigr>}{4\bigl<LE\bigr>}+\frac{\bigl<L^2\bigr>^2\bigl<E^2\bigr>}{24\bigl<LE\bigr>^3}\quad.$$}

{\it Démonstration:} On pose $\alpha=\bigl<LE\bigr>/\bigl<E^2\bigr>$ et
$\beta=\bigl<L^2\bigr>/\bigl<LE\bigr>$. Remarquons que $\beta\le\alpha$ par le
théorème de l'indice de Hodge.

Grâce au lemme 2.1, on a les estimations suivantes:
$$\begin{array}{rcl}
\D S_n=\sum_{k\ge1}h^0(X;nL-kE)&\ge&\D\sum_{k=1}^{[\beta n/2]}\Bigl(\bigl<L^2\bigr>\frac{n^2}{2}-\bigl<LE\bigr>nk+\bigl<E^2\bigr>\frac{k^2}{2}\Bigr)-O(n^2)\\
&=&\D\Bigl(\bigl<L^2\bigr>\frac{\beta}{4}-\bigl<LE\bigr>\frac{\beta^2}{8}+\bigl<E^2\bigr>\frac{\beta^3}{48}\Bigr)n^3-O(n^2)\quad.\\
\end{array}$$

D'où la minoration $\D\nu(L;E)\ge\frac{\beta}{4}+\frac{\beta^2}{24\alpha}$.
$\square$\\

Montrons maintenant le résultat principal de cette section:\\

{\bf Proposition 2.3:} {\it Soit $X$ une surface projective sur $K$. Soient
$D_1;\cdots;D_r$ des diviseurs effectifs amples sur $X$. Il existe alors des
entiers $m_1;\cdots;m_r$ tels qu'en posant $\D L=\sum_{i=1}^rm_iD_i$, on ait
$m_i\ge1$ et $\D\nu(L;D_i)>\frac{r}{4}m_i$ pour tout $i\in\{1;\cdots;r\}$.}\\

{\it Démonstration:} On pose
$\Delta=\{(t_1;\cdots;t_r)\in\mathbb{R}_+^r\ |\ t_1+\cdots+t_r=1\}$. Pour tout
$t=(t_1;\cdots;t_r)\in\Delta$, on désigne par $L_t$ le $\mathbb{R}$-diviseur
$\D L_t=\sum_{j=1}^rt_jD_j$ et on pose
$\D\phi(t)=\Bigl(\sum_{i=1}^r\frac{1}{\bigl<L_tD_i\bigr>}\Bigr)^{-1}$.

On note $f:\Delta\rightarrow\Delta$ l'application continue définie par
$\D f(t)=\Bigl(\frac{\phi(t)}{\bigl<L_tD_1\bigr>};\cdots;\frac{\phi(t)}{\bigl<L_tD_r\bigr>}\Bigr)$ pour tout $t\in\Delta$. D'après le théorème de Brouwer, $f$
admet un point fixe $x=(x_1;\cdots;x_r)$. On a alors
$\phi(x)=\bigl<L_xD_i\bigr>x_i$ pour tout
$i\in\{1;\cdots;r\}$, donc $\phi(x)r=\bigl<L_x^2\bigr>$.\\

On en déduit l'inégalité $\D\frac{\bigl<L_x^2\bigr>}{\bigl<L_xD_i\bigr>}+\frac{\bigl<L_x^2\bigr>^2\bigl<D_i^2\bigr>}{6\bigl<L_xD_i\bigr>^3}> x_ir$ pour tout
$i\in\{1;\cdots;r\}$.

On approche $x$ par un $y\in\mathbb{Q}_+^{*r}\cap\Delta$ de la forme
$\D y=\Bigl(\frac{m_1}{m};\cdots;\frac{m_r}{m}\Bigr)$ de telle sorte que
l'inégalité précédente soit encore valable avec $y$ au lieu de $x$, et on
conclut en appliquant le corollaire 2.2. $\square$\\

Terminons cette section par une définition:\\

{\bf Définition:} Soit $X$ une surface normale projective sur $K$. Soient
$D_1;\cdots;D_r$ des diviseurs effectifs non nuls sur $X$. On dit que
$D_1;\cdots;D_r$ {\bf se coupent proprement} lorsque: toute intersection de
deux quelconques d'entre eux est finie, et toute intersection de trois
quelconques d'entre eux est vide.\\

\section{Arithmétique}

Soient $K$ un corps de nombres et $S$ un ensemble fini de places de $K$. Pour
tout $v\in S$, on désigne par $K_v$ le complété de $K$ en la place $v$. On note
$O_{K;S}$ l'anneau des $S$-entiers de $K$, {\it i.e.} l'ensemble des $x\in K$
tels que $|x|_v\le1$ pour toute place finie $v\notin S$.\\

{\bf Définition:} Soit $Y$ une variété sur $K$. Un ensemble
${\cal E}\subset Y(K)$ est dit {\bf $S$-entier} sur $Y$ lorsqu'il existe un
$O_{K;S}$-schéma intègre et quasi-projectif ${\cal Y}$ de fibre générique $Y$
tel que ${\cal E}\subset{\cal Y}(O_{K;S})$.\\

On va utiliser la version suivante du théorème du sous-espace de Schmidt et
Schlickewei:\\

{\bf Proposition 3.1:} {\it Soient $X$ une variété projective sur $K$ et $L$ un
faisceau inversible très ample sur $X$. Soit ${\rm h}_L$ une hauteur de Weil
relativement à $L$. Pour chaque $v\in S$, on munit $L_v$ d'une métrique
$\|\ \|_v$ et on choisit une base $(s_{1v};\cdots;s_{qv})$ de $\Gamma(X;L)$.
Soient $c\in\mathbb{R}$ et $\varepsilon>0$. Alors l'ensemble des points
$P\in X(K)$ vérifiant
$$-\sum_{v\in S}\sum_{k=1}^q\ln\|s_{kv}(P)\|_v\ge(q+\varepsilon){\rm h}_L(P)-c\qquad(*)$$
n'est pas Zariski-dense dans $X$.}\\

{\it Démonstration:} En posant $V=\Gamma(X;L)$, on a un plongement
$X\hookrightarrow\mathbb{P}(V)$. On applique alors la reformulation de Vojta
({\it cf} théorème 2.2.4 de \cite{Voj1} p. 19) du théorème du sous-espace: il
existe une réunion finie $H$ de $K$-hyperplans de
$\mathbb{P}(V)\simeq\mathbb{P}^{q-1}_K$ telle que tout point $P\in X(K)$
vérifiant $(*)$ est dans $H\cap X$. $\square$\\

On montre ci-dessous une légère extension d'un résultat de Corvaja et Zannier
({\it cf} théorème principal de \cite{CoZa} p. 707-708); on s'inspire de leur
méthode, tout en adoptant un point de vue plus géométrique:\\

{\bf Théorème 3.2:} {\it Soit $X$ une surface normale projective sur $K$.
Soient $D_1;\cdots;D_r$ des diviseurs effectifs non nuls sur $X$ qui se coupent
proprement. Posons $Y=X-D_1\cup\cdots\cup D_r$. Soient $m_1;\cdots;m_r$ des
entiers $\ge1$. On suppose que le diviseur $\D L=\sum_{i=1}^rm_iD_i$ est ample
sur $X$ et que $\nu(L;D_i)>m_i$ pour tout $i\in\{1;\cdots;r\}$. Soit
${\cal E}\subset Y(K)$ un ensemble $S$-entier sur $Y$. Alors ${\cal E}$ n'est
pas Zariski-dense dans $Y$.}\\

{\it Démonstration:} On fixe un réel $\varepsilon>0$ tel que
$\nu(L;D_i)>(1+\varepsilon)m_i$ pour tout $i\in\{1;\cdots;r\}$, puis un entier
$b\ge1$ tel que $bL$ soit très ample et que
$$\sum_{k\ge1}h^0(X;bL-kD_i)\ge(1+\varepsilon)h^0(X;bL)m_ib\quad\mbox{pour tout}\ i\in\{1;\cdots;r\}.$$
On note $q=h^0(X;bL)$ et on choisit une hauteur de Weil ${\rm h}_{bL}$
relativement à $bL$. Pour tout diviseur effectif $E$ sur $X$, on désigne par
$1_E$ la section globale de ${\cal O}_X(E)$ qu'il définit.\\

Raisonnons par l'absurde en supposant ${\cal E}$
Zariski-dense. Il existe alors une suite $(P_n)_{n\ge0}$ d'éléments de
${\cal E}$ qui est générique, {\it i.e.} telle que pour tout fermé $Z\neq X$,
l'ensemble $\{n\in\mathbb{N}\ |\ P_n\in Z\}$ est fini.

Quitte à extraire, on peut supposer (par compacité) que pour tout $v\in S$, la
suite $(P_{nv})_{n\ge0}$ converge dans $X(K_v)$ vers un $y_v\in X(K_v)$.\\

Soit $v\in S$. On munit chaque faisceau ${\cal O}_X(D_i)_v$ d'une métrique
$\|\ \|_v$. Les diviseurs $D_1;\cdots;D_r$ se coupent proprement, donc il
existe deux indices $j_v<l_v$ tels que $y_v\notin D_i$ pour tout
$i\in\{1;\cdots;r\}-\{j_v;l_v\}$.

Le lemme 3.2 de \cite{CoZa} fournit une base
${\cal B}_v=(s_{1v};\cdots;s_{qv})$ de $\Gamma(X;bL)$ adaptée aux filtrations
$\Bigl(\Gamma(X;bL-kD_{j_v})\Bigr)_{k\ge0}$ et
$\Bigl(\Gamma(X;bL-kD_{l_v})\Bigr)_{k\ge0}$, {\it i.e.} ${\cal B}_v$ contient
une base de $\Gamma(X;bL-kD_i)$ pour tout $i\in\{j_v;l_v\}$ et tout $k\ge0$.\\

\underline{Fait}: On a la minoration suivante pour tout $n\ge0$:
$$-\sum_{k=1}^q\ln\|s_{kv}(P_n)\|_v\ge-(q+q\varepsilon)\ln\|1_{bL}(P_n)\|_v-O(1)\quad,\qquad(1)$$
où le $O(1)$ est indépendant de $n$.\\

Prouvons ce fait. Soit $s\in\Gamma(X;bL)-\{0\}$. Pour tout
$i\in\{1;\cdots;r\}$, notons $\mu_i(s)$ le plus grand entier $\mu$ tel que le
diviseur ${\rm div}(s)-\mu D_i$ soit effectif. Puisque les fermés $D_{j_v}$ et
$D_{l_v}$ n'ont pas de composante commune, le diviseur
${\rm div}(s)-\mu_{j_v}(s)D_{j_v}-\mu_{l_v}(s)D_{l_v}$ est encore effectif.
Ceci implique l'inégalité
$$-\ln\|s(P_n)\|_v\ge-\mu_{j_v}(s)\ln\|1_{D_{j_v}}(P_n)\|_v-\mu_{l_v}(s)\ln\|1_{D_{l_v}}(P_n)\|_v-O(1)\quad.$$

On écrit cette inégalité pour $s=s_{kv}$, puis on somme sur $k$. En observant
que pour $i\in\{j_v;l_v\}$, on a
$$\begin{array}{rcl}
\D\sum_{k=1}^q\mu_i(s_{kv})&=&\D\sum_{\mu\ge0}\Bigl[h^0(X;bL-\mu D_i)-h^0(X;bL-(\mu+1)D_i)\Bigr]\mu\\
&=&\D\sum_{\mu\ge1}h^0(X;bL-\mu D_i)\ge(q+q\varepsilon)m_ib\quad,\\
\end{array}$$
on trouve alors
$$-\sum_{k=1}^q\ln\|s_{kv}(P_n)\|_v\ge-(q+q\varepsilon)b\Bigl[m_{j_v}\ln\|1_{D_{j_v}}(P_n)\|_v+m_{l_v}\ln\|1_{D_{l_v}}(P_n)\|_v\Bigr]-O(1)\quad.$$

Le fait énoncé s'en déduit en remarquant que $\ln\|1_{D_i}(P_n)\|_v=O(1)$
pour tout $i\in\{1;\cdots;r\}-\{j_v;l_v\}$.\\

Maintenant, l'ensemble ${\cal E}$ est $S$-entier sur $Y$, donc pour tout
$n\ge0$, on a
$${\rm h}_{bL}(P_n)=-\sum_{v\in S}\ln\|1_{bL}(P_n)\|_v+O(1)\quad.$$

En utilisant la minoration $(1)$, on obtient (pour tout $n\ge0$)
$$-\sum_{v\in S}\sum_{k=1}^q\ln\|s_{kv}(P_n)\|_v\ge(q+q\varepsilon){\rm h}_{bL}(P_n)-O(1)\quad.$$
D'où une contradiction avec la proposition 3.1 ({\it i.e.} le théorème du
sous-espace). $\square$\\

{\it Démonstration du théorème 1.1:} Il suffit d'appliquer la proposition 2.3
(avec $r=4$) puis le théorème 3.2. $\square$\\

\ \\

{\small Pascal Autissier. I.R.M.A.R., Université de Rennes I, campus de
Beaulieu, 35042 Rennes cedex, France.

pascal.autissier@univ-rennes1.fr}

\end{document}